\documentclass[11pt]{article}

\usepackage{amsmath}
\usepackage{amsthm}
\usepackage{graphicx}
\usepackage{color}
\usepackage{amsfonts}
\usepackage{amssymb}
\usepackage{psfrag}
\usepackage{wrapfig}
\usepackage{amscd}
\usepackage{url}

\newcommand{\re}{\mathbb{R}}

\newcommand{\bdS}{\mathcal{S}}
\newcommand{\diag}{\mbox{diag}}

\newcommand{\half}{\frac{1}{2}}
\newcommand{\lmd}{\lambda}

\newcommand{\nn}{\nonumber}

\def\af{\alpha}

\newcommand{\reff}[1]{(\ref{#1})}
\newcommand{\pt}{\partial}

\newcommand{\mc}[1]{\mathcal{#1}}

\newcommand{\bdes}{\begin{description}}
\newcommand{\edes}{\end{description}}

\newcommand{\bal}{\begin{align}}
\newcommand{\eal}{\end{align}}

\newcommand{\bnum}{\begin{enumerate}}
\newcommand{\enum}{\end{enumerate}}

\newcommand{\bit}{\begin{itemize}}
\newcommand{\eit}{\end{itemize}}

\newcommand{\bea}{\begin{eqnarray}}
\newcommand{\eea}{\end{eqnarray}}
\newcommand{\be}{\begin{equation}}
\newcommand{\ee}{\end{equation}}

\newcommand{\baray}{\begin{array}}
\newcommand{\earay}{\end{array}}

\newcommand{\bsry}{\begin{subarray}}
\newcommand{\esry}{\end{subarray}}

\newcommand{\bca}{\begin{cases}}
\newcommand{\eca}{\end{cases}}

\newcommand{\bcen}{\begin{center}}
\newcommand{\ecen}{\end{center}}

\newcommand{\bbm}{\begin{bmatrix}}
\newcommand{\ebm}{\end{bmatrix}}

\newcommand{\bmx}{\begin{matrix}}
\newcommand{\emx}{\end{matrix}}

\newcommand{\bpm}{\begin{pmatrix}}
\newcommand{\epm}{\end{pmatrix}}

\newcommand{\btab}{\begin{tabular}}
\newcommand{\etab}{\end{tabular}}

\newcommand{\thmlist}{
\begin{list}{Step 1}
{\setlength{\leftmargin}{0.6 in}\setlength{\labelwidth}{0.5 in}} }
\newcommand{\alglist}{
\begin{list}{Step 1}
{\setlength{\leftmargin}{1.1 in}\setlength{\labelwidth}{1.0 in}} }
\theoremstyle{plain}
\newtheorem{theorem}{Theorem}[section]

\newtheorem{lemma}[theorem]{Lemma}

\theoremstyle{definition}
\newtheorem{example}[theorem]{Example}

\renewcommand{\subsection}[1]{
    \stepcounter{subsection}
    \settowidth{\hangindent}{\bf\thesubsection.~}
    \hangafter=1
    \bigskip\bigskip\noindent
    {\bf\hbox{\thesubsection.~}#1}\par
    \nobreak
    \medskip
}
\usepackage[top=1in,bottom=1in,left=1in,right=1in]{geometry}

\def\defg{g}

\def\hS{\hat S}

\begin{document}

\twocolumn

\title{Structured Semidefinite Representation of Some Convex Sets
\author{J.~William Helton and Jiawang Nie \\
Department of Mathematics \\
University of California, San Diego\\
helton@math.ucsd.edu,\quad
njw@math.ucsd.edu}
\date{February 10, 2008}
}

\maketitle

\begin{abstract}
Linear matrix Inequalities
(LMIs) have had a major impact on control
but formulating a problem as an LMI is an art.
Recently there is the beginnings of a  theory of which problems
are in fact expressible as LMIs.
For optimization purposes it
can also be useful to have ``lifts" which
are expressible as LMIs.
We show here that this is a much less restrictive
condition and give methods for actually constructing lifts
and their LMI representation.
\end{abstract}

\section{Introduction}

Recently, there is a lot of work \cite{Las01,ND05,NDS,Par00,ParStu}
in solving global polynomial optimization problems
by using sum of squares (SOS) methods or
semidefinite programming (SDP) relaxations.
The basic idea is to approximate
a semialgebraic set $S$ by a
collection of convex sets called {SDP relaxations}
each of which has an SDP representation.
This leads to the fundamental problem of
which sets can be represented with LMIs
or projections of LMIs.

A set $S$ is said to have an {\it LMI representation} or be {\it
LMI representable} if
 \be
  S=\{x\in \re^n: A_0+\sum_{i=1}^n A_i x_i \succeq 0\}
\ee
for some symmetric matrices $A_i$.
Here the notation $X
\succeq 0 \,(\succ 0)$ means the matrix $X$ is positive
semidefinite (definite).
Obvious necessary conditions for $S$ to be LMI representable are
that $S$ must be convex and $S$ must have the form
{\small
\be \label{eq:semialg}
S =\{x\in \re^n:\, g_1(x)\geq 0, \cdots, g_m(x)\geq 0\}
\ee
}
where $g_i(x)$ are multivariate polynomials;
such are called  basic closed semialgebraic sets.


It turns out that many convex sets are not LMI representable,
see Helton and Vinnikov \cite{HV07}.
For instance, the convex set
\[
\{x\in \re^2:\, 1-(x_1^4+x_2^4) \geq 0 \}
\]
does not admit an LMI representation.
However, the set $S$ is the projection onto $x$-space of the set
\begin{align*}
\hS:= \Big\{(x,w)\in \re^2\times \re^2:\,
\bbm 1 & x_2 \\ x_2 & w_2 \ebm \succeq 0  \\
\bbm 1+w_1 & w_2 \\ w_2
& 1-w_1 \ebm \succeq 0, \bbm 1 & x_1 \\ x_1 & w_1 \ebm \succeq 0 \Big\}
\end{align*}
in $\re^4$ which is representable by an LMI.
\medskip

More generally a set $S\subseteq \re^n$ is said to be
{\it semidefinite representable } or
{\it SDP representable} if $S$ can be described as
{\small
\begin{align} \nn
S = \Big\{x\in\re^n:\, \exists w \in \re^M \, \mbox{s.t.} \qquad \qquad  \\
\qquad \qquad A+\sum_{i=1}^n x_iB_i +\sum_{j=1}^M w_jC_j \succeq 0 \Big\}.  \nn
\end{align}
}
Here $A,B_i,C_j$ are symmetric matrices of appropriate dimensions.
Conceptually, one can think of $S$ as the
projection into $\re^n$ of a set $\hS$
in $\re^{(n+M)}$  having the LMI representation:
{\small
\begin{align}
\nn
\hS = \Big\{(x, w) \in\re^{(n+ M)}:\, \qquad \qquad \qquad \\
\qquad \qquad A+\sum_{i=1}^n x_iB_i +\sum_{j=1}^M w_jC_j \succeq 0 \Big\}.  \label{eq:SdpR}
\end{align}
The representation \reff{eq:SdpR} is called a
{\it semidefinite representation } or
{\it SDP representation} of the set $S$.
We refer to the $w_j$ as {\it auxillary variables}.

The key use of an SDP representation is
illustrated by optimizing a linear function $\ell^Tx$ over $S$.
Note that minimizing $\ell^Tx$ over $S$ is equivalent to problem
\[
\min_{(x,y) \in \hat S} \ell^Tx,
\]
which is a conventional LMI, so can be attacked by standard toolboxes.
Nesterov and Nemirovski (\cite{NN94})
in their book which introduced LMIs and  SDP
gave collections of examples of SDP representable sets
thereby leading  to:

\medskip
\noindent
{\bf Question:} \ {\it Which convex sets $S$ are the
projection of a set $\hS$ having an LMI representation?
}

\medskip

\noindent
In \S 4.3.1 of his excellent 2006 survey \cite{Nem06}, Nemirovsky commented
``this question seems to be completely open".
Now much more is known and
this paper describes both qualitative theory and SDP constructions.

%
%

Recently, Helton and Nie \cite{HN1,HN2}
proved some sufficient conditions that
guarantee the convex set $S$ is SDP representable.
For instance, one sufficient condition is called the so-called
{\it sos-convexity} or {\it sos-concavity}.
A polynomial $f(x)$ is called sos-convex
if its Hessian matrix $\nabla^2f(x) = W(x)^TW(x)$
for some matrix polynomials
($W(x)$ is not necessarily square).
A polynomial $g(x)$ is called sos-concave
if $-g(x)$ is sos-convex.
Helton and Nie \cite{HN1} proved
the following theorem:

\begin{theorem} [\cite{HN1}]
\label{thm:sos-concave}
If every $g_i(x)$ is sos-concave, then
$S$ is SDP representable.
\end{theorem}

An explicit construction of one SDP representation of $S$
when every $g_i(x)$ is sos-concave will be given in Section~\ref{sec:sos}.
For general polynomials $g_i(x)$,
the constructed SDP representation in Section~\ref{sec:sos}
is usually very big.
However, when polynomials $g_i(x)$ are sparse,
the SDP representation can be reduced to have smaller sizes.
This will be addressed in Section~\ref{sec:spar}.

There are also some  sufficient conditions other than
 sos-concavity that
guarantee the SDP representability.
For instance, when the boundary of $S$ is positively curved,
then $S$ is SDP representable.
This will be discussed in Section~\ref{sec:poscurv}.

\section{General SDP representation}
\label{sec:sos}

Suppose $S$ is a convex set given in the form
{\small
\be \nn
S =\{x\in \re^n:\, g_1(x)\geq 0, \cdots, g_m(x)\geq 0\}.
\ee
}
In this section, we assume every $g_i(x)$ is a sos-concave polynomial.
A natural SDP relaxation of $S$ is
{\small
\be \label{eq:R}
R = \big\{ x:\, \exists \, y \, \, s.t. \,\,
g(x,y) \geq 0, \, M_d(x,y) \succeq 0 \big\}.
\ee
}
Here $g(x,y)$ is a vector valued linear function and
$M_d(x,y)$ is a matrix  valued linear function
defined in what follows.
The integer $2d$ is the minimum upper bound of the degrees of $g_i(x)$.
The vector $g(x,y)$ is of the form
{\small
\[
g(x,y) = \tilde g_0 + \sum_{i=1}^n x_i \tilde g_i  +
\sum_{1< |\af| \leq 2d } y_\af \tilde g_\af
\]}
whose coefficients are such that
{\small
\[
\bbm g_1(x) \\ \vdots \\ g_m(x) \ebm =
\tilde g_0 + \sum_{i=1}^n x_i \tilde g_i  + \sum_{1< |\af| \leq 2d } x^\af \tilde g_\af.
\]}
The matrix $M_d(x,y)$ is the $d$-th order {\it moment matrix} constructed as
{\small
\be \label{eq:M}
M_d(x,y)  = A_0 + \sum_{i=1}^n x_i A_i +
 \sum_{1< |\af| \leq 2N } y_\af A_\af.
\ee}
Here the symmetric matrices $A_\af$ are such that
{\small
\[
\mathbf{m}_d(x)  \mathbf{m}_d(x)^T =
A_0 + \sum_{i=1}^n x_i A_i + \sum_{1< |\af| \leq 2d } x^\af A_\af.
\]}
The notation $\mathbf{m}_d(x)$ above denotes
the column vector of monomials with degree up to $d$, i.e.,
{\small
\[
\mathbf{m}_d(x) = \bbm 1 & x_1 & \cdots & x_1^2 & x_1x_2 & \cdots & x_n^d \ebm^T.
\]}
This construction of SDP relaxations 
of the set $S$ was proposed by Parrilo \cite{Par06} and Lasserre \cite{Las06}.
When every $g_i(x)$ is sos-concave,
Helton and Nie \cite{HN1} proved $R=S$.
This result lends itself to implementation
which
we now illustrate with two examples.
After that we improve this
 SDP construction to exploit sparsity structure
when it is present in the defining polynomials $g_i$.

\begin{example}
Consider the set $S=\{x\in\re^n: g(x)\geq 0\}$ where
\[
g(x) = 1- (x_1^4+x_2^4-x_1^2x_2^2).
\]
Direct calculation shows
\[
- \nabla^2 g(x) =
\bbm x_1 & \\ & x_2 \ebm
\underbrace{\bbm 12 & -4 \\ -4 & 12 \ebm}_{\succeq 0}
\bbm x_1 & \\ & x_2 \ebm.
\]
So $g(x)$ is sos-concave. Thus we know
$S$ can be represented by $R$ constructed in \reff{eq:R},
which in this specification becomes
\[
\baray{c}
1-y_{40}-y_{04}+y_{22} \geq 0, \\
\bbm
1      & x_1     & x_2    & y_{20} & y_{11} & y_{02} \\
x_1    & y_{20}  & y_{11} & y_{30} & y_{21} & y_{12} \\
x_2    & y_{11}  & y_{03} & y_{21} & y_{12} & y_{03} \\
y_{20} & y_{30}  & y_{21} & y_{40} & y_{31} & y_{22} \\
y_{11} & y_{21}  & y_{12} & y_{31} & y_{22} & y_{13} \\
y_{02} & y_{12}  & y_{03} & y_{22} & y_{13} & y_{04} \\
\ebm \succeq 0.
\earay
\]
The matrix above is the second order moment matrix.
In this SDP representation, there are
$12$ auxiliary variables $y_{ij}$.
\end{example}

\begin{example}
Consider the set $S=\{x\in\re^n: 1 -p(x)\geq 0\}$ where
$p$
is a homogeneous polynomial:
\[
p(x) = [x^{d}]^TB [x^{d}].
\]
Here $d>0$ is an integer and
\[
B = (b_{ij})_{1\leq i,j\leq n} \succeq 0
\]
is a symmetric matrix,
and $[x^d]$ denotes the vector
\[
[x^d] = \bbm x_1^d & x_2^d  & \cdots  & x_n^d \ebm^T.
\]
Direct calculation shows
\[
\nabla^2 p(x) =
\diag( [x^{d-1}] ) \cdot W \cdot \diag( [x^{d-1}] )
\]
where the symmetric matrix $W$ is defined to be
\[
W = d^2 B + (3d^2-2d) \diag(B).
\]
Since $B\succeq 0$, we also have $W\succeq 0$.
So $p(x)$ is sos-convex. Therefore
$S$ can be represented as
{\small
\[
\left\{x:\exists\, y,\,
1- \sum_{i,j=1}^n b_{ij} y_{d(e_i+e_j)} \geq 0,
M_d(x,y) \succeq 0 \right\}.
\]
}
Here $e_i$ denotes the $i$-th standard unit basis vector of $\re^n$
and $M_d(x,y)$ is the $d$-th order moment matrix.
\end{example}

\section{Sparse SDP representation}
\label{sec:spar}

The SDP relaxations
in \cite{Las06, HN1, HN2,Par06}
have not exploited the special structures
of polynomials
$$g_1(x),\cdots,g_m(x)$$
such as dependence of each polynomial on only a few variables
(termed  sparsity).
On the other hand,
in polynomial optimization
the sparsity structure of polynomials can be exploited
to improve the computation efficiency of
their semidefinite relaxations
\cite{KKW05,Las06spr,Nie06, ND06,Par03,WKKM06}.
In this paper we show that when
the defining polynomials for $S$ are sparse,
their structures can also be exploited to get a
``sparser" SDP representation.

%

This section
gives a structured SDP relaxation
and proves a sufficient condition
such that this structured SDP relaxation represents $S$
exactly.

Throughout this section,
we assume every polynomial $g_k(x)$ is sos-concave.
Let
$$\mc{H}=\{x\in\re^n:\, a^Tx \geq b\}\supseteq S$$
be a supporting half space
and $a^Tu=b$ for some $u\in \pt S$.
When $S$ has nonempty interior, there exists
Lagrange multipliers $\lmd_1\geq0,\cdots, \lmd_m\geq 0$ such that
{\small
\[
a = \sum_{k=1}^m \lmd_k \nabla g_k(u), \, \lmd_ig_i(u)=0, \, i=1,\cdots,m.
\]}
Helton and Nie \cite{HN1} showed that the Lagrange function
{\small
\[
a^Tx-b -\sum_{k=1}^m \lmd_kg_k(x)
\]
}
is an SOS polynomial when every polynomial $g_k(x)$ is sos-concave.

Now we suppose the polynomials $g_k(x)$ are structured such that,
for any $a,b,\lmd$, there is a decomposition such that
{\small
\[
a^Tx-b-\sum_{k=1}^m \lmd_k g_k(x) =
\phi_{\lmd}^{(1)}(x_{I_1}) + \cdots + \phi_{\lmd}^{(K)}(x_{I_K}),
\]
}
where each $\phi_{\lmd}^{(i)}(x_{I_i})$ is a polynomial in variables $x_{I_i}$.
$\{I_1,\cdots,I_K\}$ is a partition of
the index set $\{1,2,\cdots,n\}$
such that $I_i \cap I_j = \emptyset$ whenever $i\ne j$.
$x_{I_i}$ denotes the subvector of $x$ whose indices are in $I_i$.
In other words, the polynomials
\[
\phi_{\lmd}^{(1)}(x_{I_1}), \cdots , \phi_{\lmd}^{(K)}(x_{I_K})
\]
are uncoupled.

Given a polynomial $p(x)$, denote by $\mbox{supp}(p(x))$
the support of $p(x)$, i.e.,
the set of exponents of existing monomials of $p(x)$.
If $p(x)$ is SOS and has decomposition $p(x) = \sum_i q_i^2(x)$,
then it holds
\[
\mbox{supp}(q_i(x)) \subseteq \mbox{convex hull} \, \Big ( \half \mbox{supp} ( p(x) ) \Big),
\]
by Theorem~1 in Reznick \cite{Rez}.
So we define $F_i$ to be the maximum lattice set such that
\[
F_i \subseteq \mbox{convex hull} \, \Big ( \half \mbox{supp} (\phi_\lmd^{(i)} ) \Big).
\]
Now define symmetric matrices $M_\af^{j}$ as follows
{\small
\begin{align*}
 & \qquad \mathbf{m}_{F_i}(x_{I_i})  \mathbf{m}_{F_i}(x_{I_i})^T  \\
= & M_0^{(i)} + \sum_{j\in I_i} x_j M_j^{(i)} + \sum_{1< |\af| \leq 2N } x^\af M_\af^{(i)}.
\end{align*}
}
Here $\mathbf{m}_{F_i}(x_{I_i})$ denotes the vector of monomials
whose exponents lie in $F_i$.
Then define linear matrices
{\scriptsize
\be \label{eq:sparM}
M_{F_i}(x,y) = M_0^{(i)} + \sum_{j\in I_i}  x_j M_j^{(i)} + \sum_{1< |\af| \leq 2d } x^\af M_\af^{(i)}.
\ee
}

\begin{lemma}
Let $a,b,\lmd$ be the above. Then there are symmetric matrices
$W_1,\cdots, W_K \succeq 0$ such that
{\scriptsize
\[
a^Tx-b -\sum_{k=1}^m \lmd_kg_k(x)  =
\sum_{i=1}^K \mathbf{m}_{F_i}(x_{I_i})^T \cdot W_i \cdot  \mathbf{m}_{F_i}(x_{I_i}).
\]
}
\end{lemma}
\begin{proof}
By the structure assumption, we have representation
{\small
\begin{align*}
L_a(x)& :=a^Tx-b -\sum_{k=1}^m \lmd_kg_k(x) \\
& = \eta_1(x_{I_1}) + \cdots + \eta_K(x_{I_K})
\end{align*}
}
for some polynomials $\eta_1(x_{I_1}),\cdots, \eta_K(x_{I_K})$.
We know $L_a(x)$ is nonnegative polynomial
and $u$ is one global minimizer such that $L_a(u)=0$.
Let $u^{(i)}$ denote the subvector of $u$
whose coordinates correspond to the variables $x_{I_i}$.
Then $u^{(i)}$ is one global minimizer of $\eta_i(x_{I_i})$.
So we know
\[
L_a(x) =  \sum_{i=1}^k  \Big(\eta_i(x_{I_i}) - \eta_i(u^{(i)}) \Big)
\]
is SOS by Section~3 in \cite{HN1}.
In the above, fix one index $i$ and set $x^{(j)}=u^{(j)}$ for $j\ne i$,
then we can see $\eta_i(x_{I_i}) - \eta_i(u^{(i)})$ must also be SOS in $x_{I_i}$.
Furthermore, by Theorem~1 in Reznick \cite{Rez}, the polynomial $\eta_i(x_{I_i}) - \eta_i(u^{(i)})$
has the representation
{\small
\[
\eta_i(x_{I_i}) - \eta_i(u^{(i)}) =
\mathbf{m}_{F_i}(x_{I_i})^T \cdot W_i \cdot  \mathbf{m}_{F_i}(x_{I_i}),
\]
}
for some symmetric matrix $W_i\succeq 0$.
Thus the Lemma is proven.
\end{proof}

\begin{theorem}
Under the above assumptions,
the convex set $S$ has the SDP representation
{\small
\begin{align} \nn
L & = \Big\{ x \in \re^n :\, \exists \, y,\,\, s.t.\,\,  g(x,y)\geq 0, \\
& \qquad \qquad
M_{F_i}(x,y)  \succeq 0,\, i=1,\cdots,K \Big\}. \label{eq:L}
\end{align}
}
That is, $S=L$.
\end{theorem}
\begin{proof}
We have seen $S \subseteq L$. If $L\ne S$, then there must exist
some point $\hat x \in L/S$.
By the Convex Set Separation Theorem, there exists one supporting hyperplane of $S$
\[
\mc{H}=\{x\in\re^n:\, a^Tx \geq b \} \supseteq S
\]
such that $a^T u = b$ for some $u\in \pt S$ and $a^T \hat x < b$.
Consider the linear optimization problem
{\small
\begin{align*}
b = \min_{x\in\re^n} & \quad a^Tx  \\
 s.t. &\quad g_1(x)\geq 0, \cdots, g_m(x) \geq 0.
\end{align*}
}
Then $u$ is one minimizer for the above.
Let $\lmd_1\geq 0,\cdots, \lmd_m\geq 0$ be the corresponding Lagrange multipliers.
Then, by the previous lemma, we have shown
{\small
\[
\baray{c}
  \qquad a^Tx  - b - \overset{m}{\underset{i=1}{\sum}} \lmd_k g_i(x) \\
= \overset{K}{\underset{i=1}{\sum}} \mathbf{m}_{F_i}(x_{I_i})^T \cdot  W_i \cdot \mathbf{m}_{F_i}(x_{I_i})
\earay
\]
}
for some symmetric matrices $W_1,\cdots,W_K \succeq 0$.
So we have
{\scriptsize
\[
\baray{c}
b = \,\, \max  \quad \gamma \, \mbox{ s.t. } \\
 a^Tx-\gamma -\overset{m}{\underset{i=1}{\sum}} \lmd_k g_i(x) =
\overset{K}{\underset{i=1}{\sum}} \mathbf{m}_{F_i}(x_{I_i})^T \cdot  W_i \cdot \mathbf{m}_{F_i}(x_{I_i})  \\
 \lmd_1, \cdots, \lmd_m \geq 0,  W_1,\cdots,W_K \succeq 0.
\earay
\]
}
The dual of the above SOS program is
\be \nn
\min \, a^Tx \quad \, s.t.\, \,\, x \in L.
\ee
Since $\hat x \in L$, by weak duality, it holds
$ b \leq  a^T \hat x$,
which contradicts the previous assertion $a^T\hat x < b$.
\end{proof}

Now let us show some examples
for the sparse SDP representation constructed in \reff{eq:L}.

\begin{example}
Consider the convex set
{\small
\[
S = \{ x\in \re_+^n: g(x):= 1-(x_1^8+x_1^2+x_1x_2+x_2^2) \geq 0 \}.
\]
}
Obviously $g(x)$ is sos-concave.
The convex hull of $\Big( \half \mbox{supp} (g) \Big)$
contains only the following integer points:
\[
(0,0), (1,0), (2,0),(3,0),(4,0),(0,1).
\]
By the sparsity theorem, $S$ can be represented as
\[
\baray{c}
1-y_{80}-y_{20}-y_{11}-y_{02} \geq 0, \\
\bbm
1      & x_1     & x_2    & y_{20} & y_{30} & y_{40} \\
x_1    & y_{20}  & y_{11} & y_{30} & y_{40} & y_{50} \\
x_2    & y_{11}  & y_{02} & y_{21} & y_{31} & y_{41} \\
y_{20} & y_{30}  & y_{21} & y_{40} & y_{50} & y_{60} \\
y_{30} & y_{40}  & y_{31} & y_{50} & y_{60} & y_{70} \\
y_{40} & y_{50}  & y_{41} & y_{60} & y_{70} & y_{80}
\ebm \succeq 0.
\earay
\]
The matrix above is the sparse moment matrix
constructed in \reff{eq:sparM}.
There are totally $11$ auxiliary variables $y_{ij}$.
\end{example}

\begin{example}
Consider the set $S=\{x\in\re^n: 1 -p(x)\geq 0\}$ where
\[
p(x) = \sum_{i=1}^n p_i(x_i), \quad p_i(x_i) = \sum_{k=1}^{2d} \frac{x_i^k}{k !}.
\]
Obviously $p(x)$ is sos-convex, because
each univariate polynomial $p_i(x_i)$
is convex and hence sos-convex.
Thus $S$ can be represented as
{\small
\[
\baray{c}
1- \overset{n}{\underset{i=1}{\sum}} \overset{2d}{\underset{k=1}{\sum}}
 \frac{y^{(i)}_k}{k !} \geq 0 \\
H_1(x_1, y^{(i)}) \succeq 0, \cdots, H_n(x_n,y^{(n)})\succeq 0
\earay
\]
}
where $H_i(x_i,y^{(i)})$ are defined as
{\small
\[
H_i(x_i,y^{(i)}) = \bbm
1 & x_i & y_2^{(i)} &  \cdots & y_d^{(i)} \\
x_i & y_2^{(i)} & y_3^{(i)} &  \cdots & y_{d+1}^{(i)} \\
y_2^{(i)} & y_3^{(i)} &  y_4^{(i)} & \cdots & y_{d+2}^{(i)} \\
\vdots & \vdots &  \vdots & \ddots & \vdots \\
y_{d}^{(i)} & y_{d+1}^{(i)} &  y_{d+2}^{(i)} &  \cdots & y_{2d}^{(i)}
\ebm.
\]
}
The symmetric matrices $H_i(x_i,y^{(i)})$ are
sparse moment matrices constructed in \reff{eq:sparM}.
There are totally $2n(d-1)$ auxiliary variables $y_k^{(i)}$.
\end{example}

\section{Positive curvature condition}
\label{sec:poscurv}

\def\bdS{\pt S}

Section~\ref{sec:sos} and Section~\ref{sec:spar}
show the explicit construction of SDP representation
when all the defining polynomials $g_i(x)$ are sos-concave.
If some $g_i(x)$ is not sos-concave,
these constructions usually do not represent $S$.
However, there are other sufficient conditions
that guarantees $S$ is SDP representable,
which is called {\it positive curvature}.

Assume $S$ in \reff{eq:semialg} is  convex,
compact and has nonempty interior.
Denote by $\bdS$ the boundary  of $S$.
Let $Z_i=\{x:\, \defg_i(x) = 0\}$ and note
$ \bdS \subset \cup_i Z_i $.
We say   the defining functions of $S$ are {\it nondegenerate}
provided $\nabla g_i(x) \neq 0$ for all
$x \in Z_i \cap \bdS$.
%
%
The boundary of $S$ is said to have
{\it positive curvature} provided that
there exist nondegenerate defining functions $g_i$ for $S$
such that
at each
$ \, x \in \bdS \cap Z_i$
\be
\label{def:posCurvature}
-v^T\nabla^2 \defg_i(x)v > 0 , \  \ \forall\,
0\ne v \in \nabla \defg_i(x)^\perp ,
\ee
in other words, the Hessian of $\defg_i$ compressed to the
tangent space (the second fundamental form)
is negative definite.
A standard fact in geometry 
is that this does not depend on the choice of $g_i(X)$.

Obviously, necessary conditions for
$S$ to be SDP representable are that
$S$ must be convex and semialgebraic (describable
by a system of polynomial equalities or inequalities
over the real numbers).
The following, Theorem 3.3 of \cite{HN2},
goes in the direction of the converse.

\begin{theorem}
\label{thm:posCurve}
Suppose $S$ is a convex
compact  set with nonempty interior which
has nondegenerate defining polynomials
$S =\{x\in \re^n:\, \defg_1(x)\geq 0, \cdots,
 \defg_m(x)\geq 0\}$.
If the boundary $\bdS$ is positively curved,
then $S$ is SDP representable.
\end{theorem}

If $S$ is convex with  nondegenerate defining functions, then
its boundary has nonnegative curvature.
Thus the positive curvature assumption is not a huge restriction
beyond being strictly convex. The nondegeneracy assumption is
another restriction.
%

Finally comes an example where
the defining polynomial is not concave
but the boundary has positive curvature.

\begin{example}
Consider the set
{\small
\[
S = \{ x\in \re_+^n: g(x): = x_1x_2 \cdots x_n - 1 \geq 0 \}.
\]
}
We can easily see that $S$ is convex
but the defining polynomial $g(x)$ is not concave.
Note that
{\small
\[ \baray{c}
\frac{\nabla g(x)}{g(x)+1}= \bbm \frac{1}{x_1} & \frac{1}{x_2} & \cdots & \frac{1}{x_n}  \ebm^T \\
\frac{\nabla^2 g(x) }{g(x)+1}=
\bbm  0 & \frac{1}{x_1x_2} & \cdots & \frac{1}{x_1x_n} \\
 \frac{1}{x_1x_2} & 0 & \cdots & \frac{1}{x_2x_n} \\
 \vdots & \vdots & \ddots & \vdots \\
\frac{1}{x_1x_n} & \cdots & \frac{1}{x_{n-1}x_n} & 0
\ebm.
\earay\]
}
We claim that the boundary $\pt S$ has positive curvature,
%
%
which is justified by the following observation:
{\small
\[\baray{c}
-\nabla^2 g(x) + \nabla g(x) \nabla g(x)^T   \\
\succeq  (g(x)+1) \diag \left( \frac{1}{x_1^2} ,
\frac{1}{x_2^2} , \cdots, \frac{1}{x_n^2} \right)
 \succ 0,\, \forall \, x\in \pt S.
\earay \]}
Since $\pt S$ has positive curvature,
Theorem \ref{thm:posCurve}
guarantees $S$ has an SDP representation
whose construction was in Section~5 in \cite{HN1}.
\end{example}

\section{Concluding remarks}

This paper gives an explicit construction,  \reff{eq:R},
of an SDP representation for
a convex set $S$ and a sparser one \reff{eq:L}
when polynomials $g_k(x)$ are sos-concave.
There are also some other constructions
of SDP relaxations \cite{Las06,HN1,HN2} for $S$,
which are also SDP representations of $S$
when $g_k(x)$ are strictly concave
on the boundary $\pt S$ of $S$ or when
the boundary $\pt S$ has positive curvature.

In theory a hierarchy of SDP relaxations converging to $S$
within finitely many steps can be constructed
when the boundary $\pt S$ has positive curvature
(weaker than our hypothesis).
However, these refined
constructions of SDP representations
are usually more complicated than
\reff{eq:R} or \reff{eq:L}, for example,
 usually it is difficult to predict
which step of their hierarchy of relaxations represents $S$ exactly.
In contrast, the size of construction \reff{eq:R} or \reff{eq:L}
is explicit.
We refer to \cite{HN1,HN2} for more details.

%
%
%
%

\end{document}